\documentclass[reqno]{amsart}
\usepackage{amssymb,setspace}
\usepackage{ifpdf}
\ifpdf
 \usepackage[hyperindex,pagebackref]{hyperref}
\else
 \expandafter\ifx\csname dvipdfm\endcsname\relax
 \usepackage[hypertex,hyperindex,pagebackref]{hyperref}
 \else
 \usepackage[dvipdfm,hyperindex,pagebackref]{hyperref}
 \fi
\fi
\allowdisplaybreaks[4]
\numberwithin{equation}{section}
\theoremstyle{plain}
\newtheorem{thm}{Theorem}[section]
\theoremstyle{definition}
\newtheorem{dfn}{Definition}[section]
\theoremstyle{remark}
\newtheorem{rem}{Remark}[section]
\DeclareMathOperator{\td}{d\mspace{-1mu}}
\newcommand{\cmdeg}[1]{\sideset{}{_\mathrm{cm}^{#1}}\deg}
\begin{document}

\title[Complete monotonicity and completely monotonic degree]
{Complete monotonicity, completely monotonic degree, integral representations, and an inequality related to the exponential, trigamma, and modified Bessel functions}

\author[F. Qi]{Feng Qi}
\address{School of Mathematics and Informatics, Henan Polytechnic University, Jiaozuo City, Henan Province, 454010, China}
\email{\href{mailto: F. Qi <qifeng618@gmail.com>}{qifeng618@gmail.com}, \href{mailto: F. Qi <qifeng618@hotmail.com>}{qifeng618@hotmail.com}, \href{mailto: F. Qi <qifeng618@qq.com>}{qifeng618@qq.com}}
\urladdr{\url{http://qifeng618.wordpress.com}}

\author[S.-H. Wang]{Shu-Hong Wang}
\address{College of Mathematics, Inner Mongolia University for Nationalities, Tongliao City, Inner Mongolia Autonomous Region, 028043, China}
\email{\href{mailto: S.-H. Wang <shuhong7682@163.com>}{shuhong7682@163.com}}



\begin{abstract}
In the paper, the authors verify the complete monotonicity of the difference $e^{1/t}-\psi'(t)$ on $(0,\infty)$, compute the completely monotonic degree and establish integral representations of the remainder of the Laurent series expansion of $e^{1/z}$, and derive an inequality which gives a lower bound for the first order modified Bessel function of the first kind.
\end{abstract}

\thanks{This paper was typeset using \AmS-\LaTeX}

\maketitle

\section{Introduction}

In~\cite[Lemma~2]{Yang-Fan-2008-Dec-simp.tex}, the inequality
\begin{equation}\label{e-1-t-1}
\psi'(t)<e^{1/t}-1
\end{equation}
on $(0,\infty)$ was obtained and applied, where $\psi(t)$ stands for the digamma function which may be defined by the logarithmic derivative
\begin{equation*}
\psi(t)=[\ln\Gamma(t)]'=\frac{\Gamma'(t)}{\Gamma(t)}
\end{equation*}
and $\Gamma(t)$ is the classical Euler gamma function which may be defined for $\Re z>0$ by
\begin{equation*}%\label{gamma-dfn}
\Gamma(z)=\int^\infty_0t^{z-1} e^{-t}\td t.
\end{equation*}
The derivatives $\psi'(z)$ and $\psi''(z)$ of $\psi(z)$ are respectively called the tri- and tetra-gamma functions. As a whole, the derivatives $\psi^{(k)}(z)$ for $k\in\{0\}\cup\mathbb{N}$ are called the polygamma functions.
\par
The first aim of this paper is to generalize the inequality~\eqref{e-1-t-1} to complete monotonicity of a difference between both sides of~\eqref{e-1-t-1}.
\par
Our first result can be formulated as Theorem~\ref{CM-Exp-thm} below.

\begin{thm}\label{CM-Exp-thm}
The function
\begin{equation}\label{alpha-exp=psi-eq}
h(t)=e^{1/t}-\psi'(t)
\end{equation}
is completely monotonic, that is, $(-1)^{k-1}h^{(k-1)}(t)\ge0$ for $k\in\mathbb{N}$, on $(0,\infty)$ and
\begin{equation}\label{h(t)-limit=1}
\lim_{t\to\infty}h(t)=1.
\end{equation}
\end{thm}

Recently, the notion ``completely monotonic degree'' was introduced in~\cite{psi-proper-fraction-degree-two.tex}, which may be regarded as a slight but essential modification of~\cite[Definition~1.5]{Koumandos-Pedersen-09-JMAA}.

\begin{dfn}[{\cite[Definition~1]{psi-proper-fraction-degree-two.tex}}]\label{x-degree-dfn}
Let $f(t)$ be a function defined on $(0,\infty)$ and have derivatives of all orders. A number $r\in\mathbb{R}\cup\{\pm\infty\}$ is said to be the completely monotonic degree of $f(t)$ with respect to $t\in(0,\infty)$ if $t^rf(t)$ is a completely monotonic function on $(0,\infty)$ but $t^{r+\varepsilon}f(t)$ is not for any positive number $\varepsilon>0$.
\end{dfn}

For convenience, the notation
\begin{equation}
\cmdeg{t}[f(t)]
\end{equation}
was designed in~\cite[p.~9890]{psi-proper-fraction-degree-two.tex} to denote the completely monotonic degree $r$ of $f(t)$ with respect to $t\in(0,\infty)$.
\par
It was pointed out in~\cite[p.~9890]{psi-proper-fraction-degree-two.tex} that the degrees of completely monotonic functions on $(0,\infty)$ are at least zero and that if a function $f(t)$ on $(0,\infty)$ has a nonnegative completely monotonic degree then it must be a completely monotonic function on $(0,\infty)$. Equivalently speaking, a function defined on $(0,\infty)$ is completely monotonic if and only if its completely monotonic degree is not negative.
\par
The second aim of this paper is to compute the completely monotonic degree and to establish integral representations of the remainder of the Laurent series expansion of the exponential function $e^{1/z}$.
\par
Our second result may be stated as the following theorem.

\begin{thm}\label{exp=k=degree=k+1-int-thm}
For $k\in\mathbb{N}\cup\{0\}$ and $z\ne0$, let
\begin{equation}\label{exp=k=sum-eq-degree=k+1}
H_k(z)=e^{1/z}-\sum_{m=0}^k\frac{1}{m!}\frac1{z^m}.
\end{equation}
\begin{enumerate}
\item
The completely monotonic degree of $H_k(t)$ on $(0,\infty)$ meets
\begin{equation}\label{H-k(t)-degree}
\cmdeg{t}[H_k(t)]=k+1.
\end{equation}
\item
For $\Re z>0$, the function $H_k(z)$ has the integral representation
\begin{equation}\label{exp=k=degree=k+1-int}
H_k(z)=\frac1{k!(k+1)!}\int_0^\infty {}_1F_2(1;k+1,k+2;t)t^k e^{-zt}\td t,
\end{equation}
where the hypergeometric series
\begin{equation}
{}_pF_q(a_1,\dotsc,a_p;b_1,\dotsc,b_q;x)=\sum_{n=0}^\infty\frac{(a_1)_n\dotsm(a_p)_n} {(b_1)_n\dotsm(b_q)_n}\frac{x^n}{n!}
\end{equation}
for $b_i\notin\{0,-1,-2,\dotsc\}$ and the shifted factorial $(a)_0=1$ and
\begin{equation}
(a)_n=a(a+1)\dotsm(a+n-1)
\end{equation}
for $n>0$ and any real or complex number $a$.
\item
For $\Re z>0$, the function $H_k(z)$ has the integral representation
\begin{equation}\label{exp=k=degree=k+1-int-bes}
H_k(z)=\frac1{z^{k+1}}\int_0^\infty \frac{I_{k+2}\bigl(2 \sqrt{t}\,\bigr)}{t^{(k+2)/2}} e^{-zt}\td t,
\end{equation}
where
\begin{equation}\label{I=nu(z)-eq}
I_\nu(z)= \sum_{k=0}^\infty\frac1{k!\Gamma(\nu+k+1)}\biggl(\frac{z}2\biggr)^{2k+\nu}
\end{equation}
for $\nu\in\mathbb{R}$ and $z\in\mathbb{C}$ is the modified Bessel function of the first kind.
\end{enumerate}
\end{thm}

As an application of Theorems~\ref{CM-Exp-thm} and~\ref{exp=k=degree=k+1-int-thm}, the following inequality for the first order modified Bessel function of the first kind $I_1$ may be derived.

\begin{thm}\label{Bessel-2-ineq-thm}
For $t>0$, we have
\begin{equation}\label{I=1-exp-ineq}
I_1(t)>\frac{(t/2)^3}{1-e^{-(t/2)^2}}.
\end{equation}
\end{thm}

\section{Proof of Theorem~\ref{CM-Exp-thm}}

From the well known formula
\begin{equation}\label{polygamma}
\psi^{(n)}(z)=(-1)^{n+1}\int_0^{\infty}\frac{u^n}{1-e^{-u}}e^{-zu}\td u
\end{equation}
for $\Re z>0$ and $n\in\mathbb{N}$, see~\cite[p.~260, 6.4.1]{abram}, it is ready that $\lim_{t\to\infty}\psi^{(n)}(t)=0$ for $n\in\mathbb{N}$. So, the limit~\eqref{h(t)-limit=1} may be deduced immediately and, by
\begin{equation}\label{exp-frac1x-expans}
\bigl(e^{1/t}\bigr)^{(i)}=(-1)^ie^{1/t}\frac1{t^{2i}}\sum_{k=0}^{i-1}a_{i,k}t^{k}
\end{equation}
for $i\in\mathbb{N}$ and $t\ne0$, where
\begin{equation}
a_{i,k}=\binom{i}{k}\binom{i-1}{k}{k!}
\end{equation}
for $0\le k\le i-1$, in~\cite[Theorem~2.1]{exp-reciprocal-cm-IJOPCM.tex},
\begin{equation}\label{h(t)-i-der-to0}
h^{(i)}(t)=\bigl(e^{1/t}\bigr)^{(i)}-\psi^{(i+1)}(t) =(-1)^ie^{1/t}\sum_{k=0}^{i-1}\frac{a_{i,k}}{t^{2i-k}}-\psi^{(i+1)}(t)\to0
\end{equation}
for $i\in\mathbb{N}$ as $t\to\infty$.
\par
Utilizing the recurrence formula
\begin{equation}\label{abram-6.4.7}
\psi^{(n)}(z+1)=\psi^{(n)}(z)+(-1)^n\frac{n!}{z^{n+1}}
\end{equation}
in~\cite[p.~260, 6.4.7]{abram} and calculating reveal
\begin{align*}
h(t+1)-h(t)&=e^{1/(t+1)}-e^{1/t}+\psi'(t)-\psi'(t+1)\\
&=e^{1/(t+1)}-e^{1/t}+\frac1{t^2}\\
&=\frac1{t^2}+\sum_{k=0}^\infty\frac1{(k+1)!}\biggl[\frac1{(t+1)^{k+1}}-\frac1{t^{k+1}}\biggr],
\end{align*}
\begin{equation*}
[h(t+1)-h(t)]^{(i)}=(-1)^i\frac{(i+1)!}{t^{i+2}} +\sum_{k=0}^\infty\frac{(-1)^i(i+k)!}{(k+1)!k!} \biggl[\frac1{(t+1)^{i+k+1}}-\frac1{t^{i+k+1}}\biggr],
\end{equation*}
and
\begin{align*}
(-1)^i[h(t+1)-h(t)]^{(i)}&=\frac{(i+1)!}{t^{i+2}} +\sum_{k=0}^\infty\frac{(i+k)!}{(k+1)!k!} \biggl[\frac1{(t+1)^{i+k+1}}-\frac1{t^{i+k+1}}\biggr]\\
&<\frac{(i+1)!}{t^{i+2}} +\sum_{k=0}^2\frac{(i+k)!}{(k+1)!k!} \biggl[\frac1{(t+1)^{i+k+1}}-\frac1{t^{i+k+1}}\biggr]\\
&=\frac{i!}{12t^{i+3}(t+1)^{i+3}}f_i(t),
\end{align*}
where
\begin{align*}
f_i(t)&=6(i+1)t(t+1)\bigl[(t+1)^{i+2}+t^{i+2}\bigr] -12 t^2(t+1)^2\bigl[(t+1)^{i+1}-t^{i+1}\bigr]\\
&\quad-(i+1)(i+2)\bigl[(t+1)^{i+3}-t^{i+3}\bigr]\\
&=6(i+1)t(t+1)\Biggl[\sum_{\ell=0}^{i+2}\binom{i+2}{\ell}t^{\ell}+t^{i+2}\Biggr] -12 t^2(t+1)^2\sum_{\ell=0}^{i}\binom{i+1}{\ell}t^{\ell}\\
&\quad-(i+1)(i+2)\sum_{\ell=0}^{i+2}\binom{i+3}{\ell}t^{\ell}\\
&=\frac{(i-1)(i+4)(i+5)}{2}\biggl[\frac{(2-i)(i+3)}{3}t-i\biggr]t^2-i(i+1)(i+5)t\\
&\quad-\sum_{\ell=4}^{i}\biggl[(i+1)(i+2)\binom{i+3}{\ell}-6(i+1)\binom{i+3}{\ell-1} +12\binom{i+3}{\ell-2}\biggr]t^\ell\\
&=\frac{(i-1)(i+4)(i+5)}{2}\biggl[\frac{(2-i)(i+3)}{3}t-i\biggr]t^2-i(i+1)(i+5)t\\
&\quad-(1+i)(2+i)-(i+4)(i+5)\sum_{\ell=4}^{i} \frac{(i-\ell+1)(i-\ell+2)}{\ell(i-\ell+5)}\binom{i+3}{\ell-1}t^\ell
\end{align*}
and an empty sum is understood to be nil. As a result, the function $f_i(t)$ is negative and
\begin{equation}\label{(-1)=i[h(t+1)-h(t)]=(i)}
(-1)^i[h(t+1)-h(t)]^{(i)}=(-1)^i[h(t+1)]^{(i)}-(-1)^i[h(t)]^{(i)}<0
\end{equation}
for all $i\ge0$ and $t\in(0,\infty)$. Hence, by consecutive recursion and~\eqref{h(t)-i-der-to0},
\begin{multline*}
(-1)^i[h(t)]^{(i)}\ge(-1)^i[h(t+1)]^{(i)}\ge(-1)^i[h(t+2)]^{(i)}\ge \dotsm\\*
\ge(-1)^i[h(t+k)]^{(i)}\ge(-1)^i\lim_{k\to\infty}[h(t+k)]^{(i)}=0
\end{multline*}
for $i\in\mathbb{N}$ and $t\in(0,\infty)$. This implies that the function $h(t)$ is decreasing on $(0,\infty)$. Combining this monotonicity with~\eqref{h(t)-limit=1} gives $h(t)>1$ on $(0,\infty)$. In conclusion, by definition, the function $h(t)$ is completely monotonic on $(0,\infty)$. The proof of Theorem~\ref{CM-Exp-thm} is complete.

\section{Proof of Theorem~\ref{exp=k=degree=k+1-int-thm}}

It is general knowledge that the exponential function $e^{1/z}$ for $z\in\mathbb{C}$ with $z\ne0$ can be expanded into the Laurent series
\begin{equation}\label{exp-reciproc-series}
e^{1/z}=\sum_{m=0}^\infty\frac1{m!}\frac1{z^m},\quad z\ne0.
\end{equation}
Therefore, it is clear that
\begin{equation}\label{exp-residue-term}
H_k(z)=\sum_{m=k+1}^\infty\frac{1}{m!}\frac1{z^m},\quad z\ne0
\end{equation}
and $x^{k+1}H_k(x)$ is completely monotonic on $(0,\infty)$. That is,
\begin{equation}\label{H-k(t)-degree>k+1}
\cmdeg{t}[H_k(t)]\ge k+1.
\end{equation}
Since, for any $\varepsilon>0$, the function
\begin{equation*}
x^{k+1+\varepsilon}H_k(x)=x^\varepsilon\sum_{m=0}^\infty\frac{1}{(m+k+1)!}\frac1{x^m}
\end{equation*}
tends to $\infty$ as $x\to\infty$, we see that for any $\varepsilon>0$ the function $x^{k+1+\varepsilon}H_k(x)$ is not completely monotonic on $(0,\infty)$. That is,
\begin{equation}\label{H-k(t)-degree<k+1}
\cmdeg{t}[H_k(t)]\le k+1.
\end{equation}
Combining~\eqref{H-k(t)-degree>k+1} and~\eqref{H-k(t)-degree<k+1} leads to~\eqref{H-k(t)-degree}.
\par
For $\Re z>0$ and $\Re k>0$, it was listed in~\cite[p.~255, 6.1.1]{abram} that
\begin{equation*}%\label{Gamma(z)=k-z-int}
\Gamma(z)=k^z\int_0^\infty t^{z-1}e^{-kt}\td t.
\end{equation*}
This formula can be rearranged as
\begin{equation}\label{Gamma(z)=k-z-int-rearr}
\frac1{z^w}=\frac1{\Gamma(w)}\int_0^\infty t^{w-1}e^{-zt}\td t
\end{equation}
for $\Re z>0$ and $\Re w>0$. Substituting the formula~\eqref{Gamma(z)=k-z-int-rearr} into~\eqref{exp-residue-term} yields
\begin{align*}
H_k(z)&=\int_0^\infty \Biggl[\sum_{m=k+1}^\infty\frac{1}{m!}\frac1{\Gamma(m)}t^{m-1}\Biggr]e^{-zt}\td t\\
&=\frac1{k!(k+1)!}\int_0^\infty {}_1F_2(1;k+1,k+2;t)t^k e^{-zt}\td t.
\end{align*}
The integral representation~\eqref{exp=k=degree=k+1-int} follows.
\par
The function $H_k(z)$ can be rewritten as
\begin{align*}
H_k(z)&=\frac1{z^{k+1}}\sum_{m=0}^\infty\frac{1}{(k+1+m)!}\frac1{z^m}\\
&=\frac1{z^{k+1}}\int_0^\infty \Biggl[\sum_{m=0}^\infty\frac{1}{(k+1+m)!} \frac1{\Gamma(m)}t^{m-1}\Biggr] e^{-zt}\td t\\
&=\frac1{z^{k+1}}\int_0^\infty \frac{I_{k+2}\bigl(2 \sqrt{t}\,\bigr)}{t^{(k+2)/2}} e^{-zt}\td t.
\end{align*}
The integral representation~\eqref{exp=k=degree=k+1-int-bes} follows.
Theorem~\ref{exp=k=degree=k+1-int-thm} is thus proved.

\section{Proof of Theorem~\ref{Bessel-2-ineq-thm}}

When $k=0$, the integral representations~\eqref{exp=k=degree=k+1-int} and~\eqref{exp=k=degree=k+1-int-bes} become
\begin{equation}\label{open-answer-1}
e^{1/z}=1+\int_0^\infty \frac{I_1\bigl(2\sqrt{t}\,\bigr)}{\sqrt{t}\,} e^{-zt}\td t
\end{equation}
and
\begin{equation}\label{open-answer-2}
e^{1/z}=1+\frac1{z}\int_0^\infty \frac{I_{2}\bigl(2 \sqrt{t}\,\bigr)}{t} e^{-zt}\td t
\end{equation}
for $\Re z>0$. Hence, by~\eqref{polygamma} for $n=1$, the function $h(z)$ defined by~\eqref{alpha-exp=psi-eq} has the following integral representation
\begin{equation}\label{h(t)-int-rep-eq}
h(z)=1+\int_0^\infty \biggl[\frac{I_1\bigl(2\sqrt{u}\,\bigr)}{\sqrt{u}\,} -\frac{u}{1-e^{-u}}\biggr]e^{-zu}\td u
\end{equation}
and
\begin{equation}\label{h(t)-int-rep-eq-der}
(-1)^kh^{(k)}(t)=\int_0^\infty \biggl[\frac{I_1\bigl(2\sqrt{u}\,\bigr)}{\sqrt{u}\,} -\frac{u}{1-e^{-u}}\biggr]u^ke^{-tu}\td u
\end{equation}
for $k\ge1$ are completely monotonic on $(0,\infty)$.
The famous Hausdorff-Bernstein-Widder Theorem~\cite[p.~161, Theorem~12b]{widder} states that a necessary and sufficient condition that $f(x)$ should be completely monotonic for $0<x<\infty$ is that
\begin{equation} \label{berstein-1}
f(x)=\int_0^\infty e^{-xt}\td\alpha(t),
\end{equation}
where $\alpha(t)$ is non-decreasing and the integral converges for $0<x<\infty$. Consequently, the function in the bracket of~\eqref{h(t)-int-rep-eq-der} is not less than zero, that is,
\begin{equation*}
\frac{I_1\bigl(2\sqrt{u}\,\bigr)}{\sqrt{u}\,} \ge\frac{u}{1-e^{-u}}
\end{equation*}
in which replacing $2\sqrt{u}\,$ by $t$ yields the inequality~\eqref{I=1-exp-ineq}. The proof of Theorem~\ref{Bessel-2-ineq-thm} is complete.

\begin{rem}
The integral representations~\eqref{open-answer-1} and~\eqref{open-answer-2} supply answers to an open problem posed in~\cite[p.~127, Section~4]{exp-reciprocal-cm-IJOPCM.tex}.
\end{rem}

\begin{rem}
This paper is a corrected and extended version of the preprint~\cite{simp-exp-degree.tex}.
\end{rem}

\end{document}